\documentclass[11pt]{amsart}

\usepackage{fullpage, parskip}
\usepackage{amsmath, amsthm, amsfonts}
\usepackage{graphicx}
\usepackage{pythonhighlight}

\newtheorem{theorem}{Theorem}[section]
\newtheorem{corollary}[theorem]{Corollary}
\newtheorem{observation}[theorem]{Observation}
\newtheorem{lemma}[theorem]{Lemma}

\newcommand{\ovr}[1]{\overrightarrow{#1}}

\title{Quasi-Transitive Mixed Graphs and Undirected Squares of Oriented Graphs}
\author{Christopher Duffy}
\address{School of Mathematics and Statistics, University of Melbourne, Australia}
\begin{document}

\begin{abstract}
We consider the problem of classifying those graphs that arise as an undirected square of an oriented graph by generalising the notion of quasi-transitive directed graphs to mixed graphs.
We fully classify those graphs of maximum degree three and graphs of girth at least four that arise an undirected square of an oriented graph.
In contrast to the recognition problem for graphs that admit a quasi-transitive orientation, we find it is NP-complete to decide if a graph admits a partial orientation as a quasi-transitive mixed graph.
We prove the problem is Polynomial when restricted to inputs of maximum degree three, but remains NP-complete when restricted to inputs with maximum degree at least five.
Our proof further implies that for  fixed $k \geq 3$, it is NP-complete to decide if a graph arises as an undirected square of an orientation of a graph with $\Delta = k$.
\end{abstract}

\maketitle

\section{Introduction and Background}

Let $\Gamma$ be a graph.
Recall that $\Gamma$ is a \emph{graph square} when there exists a graph $\Sigma$ such that $\Gamma = \Sigma^2$, where $\Sigma^2$ is the graph formed from $\Sigma$ by adding an edge between any pair of vertices at distance two.
We extend this notion to oriented graphs, that is, those graphs that arise from simple graphs by assigning each edge an orientation as an arc.

We say a graph $\Gamma$ is an \emph{oriented graph square} when there exists an oriented graph $\ovr{G}$ such that $\Gamma = U(\ovr{G}^2)$, where $U(\cdot)$ denotes the simple graph underlying a mixed graph, and $\ovr{G}^2$ is the mixed graph formed from an oriented graph $\ovr{G}$ by adding an edge between any pair of vertices at directed distance two.
Similar to the problem of finding square roots of graphs, given $\Gamma$ it is not at all clear how one may recover $\ovr{G}$, nor when  $\ovr{G}$ is unique.
Further, given an arbitrary graph $\Sigma$ it is not clear how one can determine if $\Sigma$ is an oriented graph square.
We consider these problems herein.

To study oriented graph squares, we re-frame the problem as a graph orientation problem.
Recall that an oriented graph is \emph{quasi-transitive} when it contains no induced directed path of length two (i.e., no $2$-dipath).

\begin{theorem}\cite{GH62}
	A graph $\Gamma$ admits a quasi-transitive orientation if and only if $\Gamma$ is a comparability graph.
\end{theorem}

\begin{figure}[h]
	\includegraphics[scale=0.4]{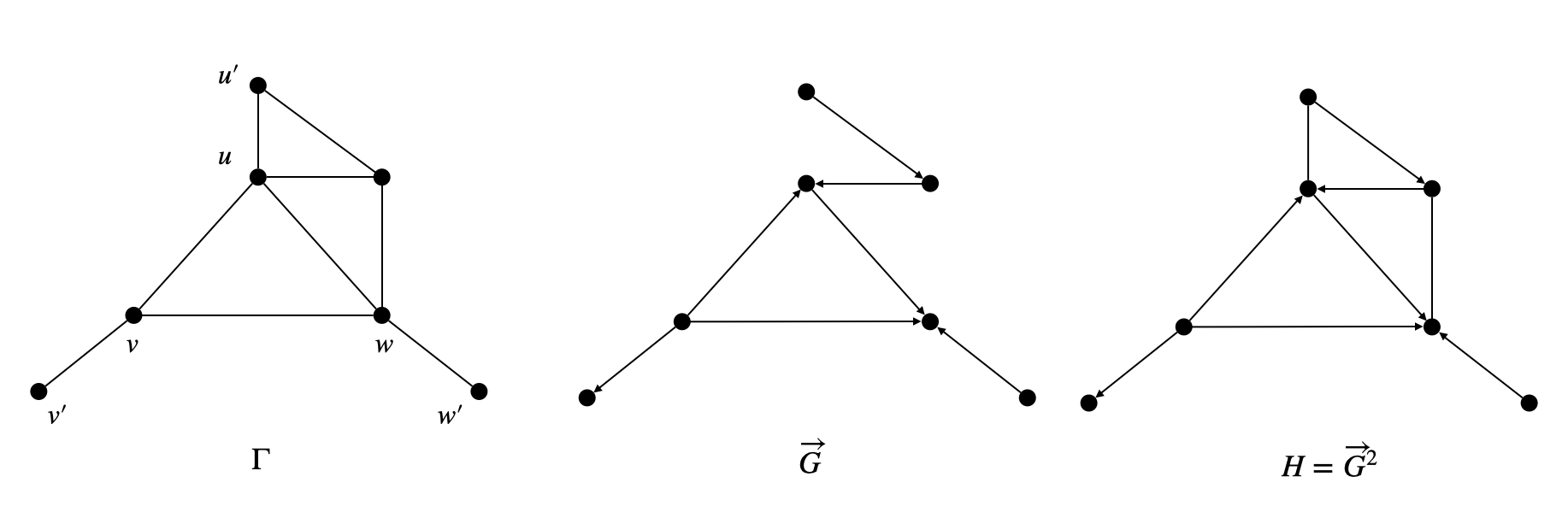}
	\caption{The oriented graph square $\Gamma$ resulting from the oriented graph $\protect\overrightarrow{G}$. $H$ is the corresponding mixed graph.}
	\label{fig:exampleFig}
\end{figure}

By definition, if $\ovr{G}$ is a quasi-transitive oriented graph, then $\ovr{G}^2 = \ovr{G}$ and so $U(\ovr{G}) = U(\ovr{G}^2)$.
This implies that for every comparability graph $\Gamma$ there exists an oriented graph $\ovr{G}$ such that $\Gamma = U(\ovr{G}^2)$.
However, there are non-comparability graphs for which there exists an oriented graph $\ovr{G}$ such that $\Gamma = U(\ovr{G}^2)$.
Consider the graphs given in Figure \ref{fig:exampleFig}.

The subgraph of $\Gamma$ induced by $\{u,v,w,u^\prime,v^\prime,w^\prime\}$ is a forbidden subgraph for the family of comparability graphs \cite{G88}.
Therefore $\Gamma$ does not admit an orientation as a quasi-transitive oriented graph.
However, $\Gamma = U(\ovr{G}^2)$.
We generalise the notation of quasi-transitive oriented graph to mixed graphs as follows.

Let $H$ be a mixed graph. 
We say $H$ is \emph{quasi-transitive} when 
\begin{enumerate}
	\item $H$ has no induced $2$-dipath; and 
	\item for every edge $uv$ in $H$ there exists $w \in V(H)$ such that $uwv$ or $vwu$ is a $2$-dipath.
\end{enumerate}
Looking back at our example, we see that $H$ is a quasi-transitive mixed graph.

Let $H$ be a quasi-transitive mixed graph and let $\Gamma = U(H)$.
Let $\ovr{G}$ be the oriented graph formed from $H$ by removing all its edges.
By the definition of quasi-transitive mixed graph we have $H = \ovr{G}^2$ and thus $\Gamma = U(\ovr{G}^2)$.
And so we arrive at the following observation.

\begin{observation}\label{obs:equiv}
	Let $\Gamma$ be a graph.
	There exists an oriented graph $\ovr{G}$ such that $\Gamma = U(\ovr{G}^2)$ if and only if $\Gamma$ admits a partial orientation as a quasi-transitive mixed graph.
\end{observation}

Observation \ref{obs:equiv} allows us to treat equivalently graphs that admit a partial orientation as a quasi-transitive mixed graphs and graphs arising as undirected squares of oriented graphs and thus unifies two previously unrelated areas: the study of (di)graph square roots and the study of quasi-transitivity in digraphs.
Throughout the remainder of this work we treat these two concepts interchangeably, opting for the one which eases the path of the reader through a particular argument or idea or which allows a direct comparison to previous work.

Quasi-transitive oriented graphs were first studied by Ghouil\`a-Houri \cite{GH62} and Gallai \cite{G67} in the 1960s.
Since then a large body of research has centred around these oriented graphs due to their connection to variety of other concepts in the study of directed graphs such as local tournaments \cite{H95} and semi-complete digraphs \cite{B92}.
Notable among these results is a proof that one can decide in polynomial time if a graph admits a quasi-transitive orientation \cite{G88}.

The study of graph squares and related concepts are a mainstay of research in structural graph theory.
Graph squares were fully classified in 1967 by Mukhopadhyay \cite{M67}.
Classifying the complexity of the recognition problem for graph squares remained open until 1994 when Motwani and Sudan proved the problem is NP-complete  \cite{M94}. 
However, there is a full characterisation of tree squares (that is, graphs that arise as a square of a tree), which leads to a polynomial time algorithm for their recognition \cite{R60,L95}.
More recently the Farzad et al. presented a full dichotomy theorem for the problem of recognising if a graph is a graph square of a graph with fixed girth.
It is NP-complete to decide if a graph is a graph square of graph with girth at most $5$ and Polynomial to decide if a graph is a graph square of a graph with girth $g \geq 6$ \cite{A09,F12}.
Despite the abundance of research into graph squares, relatively little work has been devoted to the study of squares of directed graphs.
In \cite{Gl68} Geller gives a classification of digraph squares from which Mukhopadhyay's result follows as a corollary.

The remainder of this work proceeds as follows.
In the subsequent section we develop necessary lemmas and preliminary results required for the main results in the subsequent sections.
In doing so we immediately arrive at a classification of oriented graph squares of girth four.
In Section \ref{sec:cubic} we classify graphs with maximum degree three that admit a partial orientation as a quasi-transitive mixed graphs and thus classify oriented graph  squares of maximum degree three.
In Section \ref{sec:complexity} we study the computational complexity of the recognition problem for undirected squares and graphs admitting a partial orientation as a quasi-transitive mixed graph.
We prove the problem of deciding if a graph is an oriented graph square is NP-complete.
Finally in Section \ref{sec:discussion} we provide additional discussion, which relates the problem under study herein to the study of the oriented chromatic number as well as highlights potential areas of future study.

At times throughout we will be interested in the graph formed from a subset of the edges of a graph.
Let $\Gamma$ be a graph and let $X \subseteq E(\Gamma)$.
The graph \emph{formed from $X$} is the graph whose vertex set is those vertices that appear as an endpoint of an edge in $X$ and whose edge set is $X$.
Acknowledging the slight abuse of notation, we refer to this graph as $\Gamma[X]$.
We extend this convention to subsets of arcs and edges in mixed and oriented graphs.

Throughout this work we will be interested in various properties and parameters of graphs, mixed graphs and oriented graphs.
For ease of communication, we may refer to an undirected graph property or parameter of a mixed or oriented graph to mean the property or parameter of the underlying simple graph.
For example, when we say that a vertex has degree $k$ in $\ovr{G}$ we mean it has degree $k$ in $U(\ovr{G})$.
To avoid needing to consider edge cases in which a mixed graph has no edges, we consider an oriented graph to be a partial orientation of graph, that is, a mixed graph.
We use $E(G)$ and $A(G)$ to respectively refer to the edge set and the arc set of a mixed graph.
And so we may consider an oriented graph as a mixed graph $G$ with $E(G) = \emptyset$.
We say that a vertex $v$ is a source (respectively, a sink) in a mixed graph $G$ when $v$ is a source (respectively sink) in $G[A(G)]$.
For other graph theoretic notation not defined herein we refer the reader to \cite{bondy}.

\section{Preliminaries}\label{sec:prelim}
The classification of quasi-transitive oriented graphs as arising from comparability graphs gives a forbidden subgraph classification for the family of graphs that admit such an orientation \cite{G88}.
We begin by proving no such classification exists for quasi-transitive mixed graphs.

\begin{theorem}
	For every graph $\Gamma$ there exists a graph $\Gamma^\prime$ such that $\Gamma^\prime$ admits a partial orientation as a quasi-transitive mixed graph and $\Gamma$ is an induced subgraph of $\Gamma^\prime$.
\end{theorem}

\begin{proof}
	Let $\Gamma$ be a graph.
	Let $\ovr{G}$ be the oriented graph formed from $\Gamma$ by arbitrarily orienting each edge into an arc, and then bisecting each arc to be a $2$-dipath.
	Let $\Gamma^\prime = U(\ovr{G}^2)$.
	By construction $\Gamma^\prime[V(\Gamma)] = \Gamma$.
	By Observation \ref{obs:equiv}, $\Gamma^\prime$ admits a partial orientation as a quasi-transitive mixed graph.
\end{proof}

Though in general no forbidden subgraph characterisation exists, we do find such characterisations when we consider various restrictions on $\Gamma$.
Below we prove such a characterisation exists when $\Gamma$ has girth at least four.
In the subsequent section we prove such a characterisation exists when $\Gamma$ has maximum degree three.

Throughout the remainder of this work we will apply the following lemmas, which save us repetitious arguments in subsequent arguments. 

\begin{lemma}\label{lem:noK3}
	Let $\Gamma$ be an oriented graph square.
	Let $e \in E(\Gamma)$ be contained in no copy of $K_3$ in $\Gamma$.
	In every partial orientation $G$ of $\Gamma$ as a mixed quasi-transitive graph,
	\begin{enumerate}
		\item $e$  is oriented as an arc;
		\item the ends of $e$ are respectively a source or sink; and
		\item the undirected graph formed the set of edges contained in no copy of $K_3$ in $\Gamma$ contains no odd cycle.
	\end{enumerate}
\end{lemma}

\begin{proof}
	Let $\Gamma$ be a graph arising as the square of an oriented graph.
	Let $G$ be a partial orientation of $\Gamma$ as a quasi-transitive mixed graph.
	Let $X$ be set of edges not contained in a copy of $K_3$ in $\Gamma$.
	Consider $e \in X$ with $e = xy$.
	Since there is no vertex $z$ such that $zx,yz \in E(\Gamma)$, there is no $2$-dipath in $G$ with ends at $x$ and $y$.
	Therefore $e$ is oriented, in some direction, as an arc in $G$.
	
	Consider a pair of distinct edges $e_1, e_2 \in X$ with a common end-point.
	Let $e_1 = uv$ and $e_2  = vw$
	By construction of $X$, $wu \notin E(\Gamma)$.
	Since $G$ is quasi-transitive mixed graph, $uvw$ is not an induced $2$-dipath in $G$.
	By construction, $uvw$ is not an induced $2$-dipath in $G[X]$.
	
	Therefore $G[X]$ contains no induced $2$-dipath.
	And so $G[X]$ is a quasi-transitive oriented graph.
	Therefore $U(G[X])$ is a comparability graph.
	Recall that every comparability graph is perfect (see \cite{G04})
	By construction, $U(G[X])$ contains no copy of $K_3$.
	Since $U(G[X])$ is perfect and contains no copy of $K_3$, necessarily $\chi(\Gamma) \leq 2$.
	And so it follows $U(G[X])$ is bipartite.
\end{proof}

Using this lemma we fully classify those graphs with girth at least four that admit a partial orientation as a quasi-transitive mixed graph.

\begin{theorem}\label{thm:bigGirth}
	Let $\Gamma$ be an graph with girth $g \geq 4$.
	The graph $\Gamma$ admits a partial orientation as a quasi-transitive mixed graph if and only if $\Gamma$ has no odd cycle.
\end{theorem}

\begin{proof}
	Let $\Gamma$ be an graph with girth $g \geq 4$.
	If $\Gamma$ admits a partial orientation as a quasi-transitive mixed graph, then by Lemma \ref{lem:noK3}, $\Gamma$ is bipartite and thus has no odd cycle.
	
	Let $\Gamma$ be bipartite with partition $\{V_1, V_2\}$.
	We find a partial orientation of $\Gamma$ as a quasi-transitive oriented graph $G$ by orienting all arcs to have their tail in $V_1$ and their head in $V_2$.
	In such an orientation every vertex is either a source vertex or a sink vertex.
	The mixed graph $G$ has no induced $2$-dipath and no edges.
	Therefore $G$ is a partial orientation as a quasi-transitive mixed graph.
\end{proof}

Let $\Gamma$ be a graph that admits a partial orientation as a quasi-transitive mixed graph and consider the existence of a cut vertex $v$.
Let $\{V_1,V_2\}$ be a partition of $V(\Gamma)\setminus \{v\}$ such that $v$ has a neighbour in both $V_1$ and $V_2$ and no vertex of $V_1$ is adjacent to a vertex of $V_2$.
Let $G$ be a partial orientation of $\Gamma$ as a quasi-transitive mixed graph.

By construction, $v$ cannot be the centre of a $2$-dipath with an end in $V_1$ and an end in $V_2$.
Furthermore, if $v$ is the centre of a $2$-dipath with both ends in $V_1$, then no edge between $v$ and a vertex of $V_2$ can be oriented as an arc in $G$.
However this cannot be, as any vertex in an quasi-transitive mixed graph that is incident with an edge must also be an end-point of a $2$-dipath whose end is at the other end of the edge.
And so $v$ cannot be the centre of a $2$-dipath with both ends in $V_1$.
And so we conclude that $v$ is a source or a sink in $G$.
(Recall that a source/sink vertex in a mixed graph is one for which all incident arcs are oriented toward/away from the vertex.
A source/sink vertex in a mixed graph may be incident with edges)

We extend this argument to vertex cuts consisting of independent sets.
Let $\Gamma$ be a graph and let $I \subset V(\Gamma)$ be an independent set.
We say $I$ is an \emph{independent vertex cut} when there exists a partition $\{V_1,V_2,I\}$ of $V(\Gamma)$ such that $\Gamma-I$ is not connected, every vertex in $I$ is adjacent to at least one vertex in $V_1$ and one in $V_2$ and no vertex of $V_1$ is adjacent to a vertex of $V_2$.
We can extend partial orientations of $\Gamma[V_1 \cup I]$ and $\Gamma[V_2\cup I]$ as quasi-transitive mixed graphs to one to $\Gamma$ provided each vertex in the independent set is respectively a source or a sink in both the partial orientation of  $\Gamma[V_1 \cup I]$ and of $\Gamma[V_2\cup I]$.

\begin{lemma}[The Vertex Cut Lemma]
	Let $\Gamma$ be a graph and let $I \subset V(\Gamma)$ be an independent vertex cut with $I = \{v_1,v_2,\dots, v_k\}$.
	The graph $\Gamma$ admits a partial orientation as a quasi-transitive mixed graph if and only if there exists partial orientations of $\Gamma[V_1 \cup I]$ and $\Gamma[V_2\cup I]$ as quasi-transitive mixed graphs $G_1$ and $G_2$ in which for each $1 \leq i \leq k$, $v_i$ is a source (respectively, a sink) in both $G_1$ and $G_2$.
\end{lemma}

\begin{proof}
	Let $\Gamma$ be a graph and let $I \subset V(\Gamma)$ be an independent vertex cut with $I = \{v_1,v_2,\dots, v_k\}$.
	
	Let $G$ be a partial orientation of $\Gamma$ as a quasi-transitive mixed graph.
	Consider $v \in I$
	We claim $v$ is a source or a sink.
	We prove $v$ is not the centre vertex of a $2$-dipath in $G$.
	Let $u $and $u^\prime$ be incident with $v$ and  such that $u$ and $u^\prime$ are, without loss of generality, contained in $V_1$.
	
	By hypothesis, there exists $w \in V_2$ such that $w$ is adjacent to $v$.
	Further we may assume $vw$ is an arc (in some direction) in $G$, as otherwise there exists a $2$-dipath $ww^\prime v$ such that $w^\prime$ is in the same component as $w$ in $G-v_i$.
	In this case we may swap the roles of $w$ and $w^\prime$.
	
	If $uvu^\prime$ is a $2$-dipath (in some direction), then one of $u$ and $u^\prime$ must be adjacent to $w$, which contradicts our choice of $V_1$ and $V_2$.
	Therefore $v$ is not the centre of a $2$-dipath whose ends are in $V_1$.
	Therefore $v$ is a source or a sink vertex in $G[V_1 \cup I]$.
	
	By choice of $v$, it then follows that for all $v \in I$, vertex $v$ is a source or a sink vertex in $G[V_1 \cup I]$.
	Therefore there is no $2$-dipath in $G$ with an end in $V_1$ and an end in $V_2$.
	And so $G[V_1 \cup I]$ is partial orientation of $\Gamma[V_1 \cup I]$ as a quasi-transitive mixed graph.
	Similarly, $G[V_2 \cup I]$ is partial orientation of $\Gamma[V_2\cup I]$ as a quasi-transitive mixed graph.
	
	Assume now there exists partial orientations of $\Gamma[V_1 \cup I]$ and $\Gamma[V_2\cup I]$ as quasi-transitive mixed graphs $G_1$ and $G_2$ in which for each $1 \leq i \leq k$, $v_i$ is a source (respectively a sink) in both $G_1$ and $G_2$.
	Let $G$ be the mixed graph formed from $G_1$ and $G_2$ by identifying  for all $v \in I$ vertex $v$ in the copy of $G_1$ with $v$ in the copy of $G_2$.
	Since all identified vertices are either source vertices or sink vertices, this process creates nor removes $2$-dipaths.
	Similarly, this process neither creates nor removes edges.
	Since this constructions preserves the existence of all $2$-dipaths and edges, $G$ is a partial orientation of $\Gamma$ as a quasi-transitive mixed graph.		
\end{proof}

for $|I| = 1$ we arrive at the back to special case where $\Gamma$ has a cut vertex.
Thus the Vertex Cut Lemma allows us to classify cut vertices in a graph arising as an oriented graph square.

\begin{theorem}\label{thm:classifyCutVx}
	Let $\Gamma$ be a graph arising as an oriented graph square.
	If $v \in V(\Gamma)$ is a cut vertex, then $v$ is a source or a sink in every oriented graph $\ovr{G}$ such that $\Gamma = U(\ovr{G}^2)$.
\end{theorem}

\section{Quasi-Transitive Mixed Graphs with Maximum Degree Three}\label{sec:cubic}
We turn now to the problem of classifying graphs with maximum degree three that admit a partial orientation as a quasi-transitive mixed graph.
By Observation \ref{obs:equiv}, these are exactly the family of  oriented graph squares with maximum degree.
We begin by identifying a reduction that can be applied to graphs with maximum degree three that preserves the property of (not) admitting a partial orientation as a quasi-transitive mixed graph.

\begin{figure}[h]
	\includegraphics[scale=0.4]{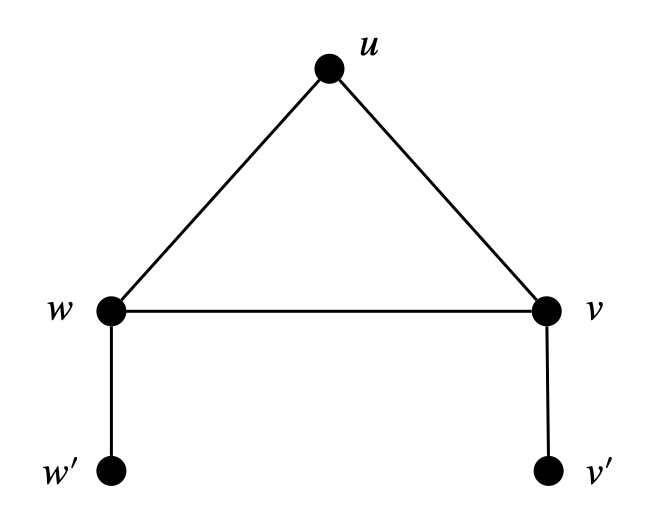}
	\caption{A removable vertex, $u$, in a graph with maximum degree three}
	\label{fig:removeableFig}
\end{figure}

Let $\Gamma$ be a graph with maximum degree three and let $u \in V(\Gamma)$.
We say $u$ is \emph{removable} when it has degree two, and its neighbours are adjacent, have degree three and have only one common neighbour (see Figure \ref{fig:removeableFig}).

\begin{lemma}\label{lem:remove1}
	Let $\Gamma$ be a graph with maximum degree three and let $u$ be a removable vertex of $\Gamma$.
	The graph $\Gamma$ admits a partial orientation as a quasi-transitive mixed graph if and only if $\Gamma - u$ admits a partial orientation as a quasi-transitive mixed graph.
\end{lemma}

\begin{proof}
	Let $\Gamma$ be a graph with maximum degree three and let $u$ be a removable vertex in $\Gamma$ with neighbours $v$ and $w$.
	Let $v^\prime$ and $w^\prime$ respectively be the other neighbour of $v$ and $w$.
	Since $\Gamma$ has maximum degree three, edges $ww^\prime$ and $vv^\prime$ are contained in no copy of $K_3$.
	
	Let $G$ be a partial orientation of $\Gamma$ as a quasi-transitive mixed graph.
	By Lemma \ref{lem:noK3}, $u$ and $v$ are each a source or a sink.
	Therefore neither $v$ nor $w$ is the centre vertex of a $2$-dipath.
	Since $v$ is not the centre of a $2$-dipath, the edge $uw$ is oriented in some direction in $G$.
	Similarly, the edge $uv$ is oriented in some direction in $G$
	
	If both of these arcs are oriented to have their head at $u$, then both $v$ and $w$ are source vertices.
	However, in this case notice that the edge $vw$ must be an edge in $G$, but there is no $2$-dipath that between $v$ and $w$.
	
	Therefore, without loss of generality, $vuw$ is a $2$-dipath in $G$, which then implies $v$ is a source and $w$ is a sink.
	If $vw$ is an arc in $G$, it must have its head at $w$.
	Furthermore, if $vw$ is an edge in $G$ it can be replaced by the arc $vw$.
	And so we may assume $vw$ is an arc in $G$.
	Removing $u$ from $G$ yields a partial orientation of $\Gamma-u$ as quasi-transitive mixed graph. 
	
	Let $H$ be a partial orientation of $\Gamma-u$ as a quasi-transitive mixed graph.
	Since none of $vv^\prime, vw$ or $ww^\prime$ are contained in a copy of $K_3$ in $\Gamma-u$, by Lemma \ref{lem:noK3}, each of these edges are oriented as arcs in $H$.
	The path $v^\prime vww^\prime$ must be oriented to have no induced $2$-dipath.
	Thus it is oriented as an alternating path.
	Without loss of generality, we assume $v^\prime v, wv, ww^\prime \in A(H)$.
	We extend the partial orientation of $\Gamma-u$ to one of $\Gamma$ by orienting the edge $uv$ to have its head at $v$ and the edge $uw$ to have its head at $u$.
\end{proof}

This lemma implies that removing a removable vertex from a graph does not affect whether it admits an partial orientation as a  quasi-transitive oriented graph.
Subsequently we prove that removing the set of removable vertices neither creates new removable vertices nor affects whether the graph admits an partial orientation as a  quasi-transitive oriented graph.

\begin{lemma}\label{lem:removePreserve}
	Let $\Gamma$ be a graph with maximum degree three and let $R$ be the set of removable vertices $\Gamma$.
	For every $u \in R$ the set of removable vertices in $\Gamma-u$ is $R\setminus u$.
\end{lemma}

\begin{proof}
		Let $\Gamma$ be a graph.
		Let $u$ and $v$ be distinct removable vertices in $\Gamma$.
		Let $u_1$ and $u_2$ be the neighbours of $u$ in $\Gamma$.
		Similarly define $v_1$ and $v_2$.
		Since $\deg(u) = 2$ and $\deg(v_1)=\deg(v_2) =3$ we have $u \neq v_1,v_2$.
		
		If $u_1 = v_1$, then since $u\neq v$, the neighbours of $u_1$ are $u,v,u_2$ and $v_2$.
		Since $u_1$ has degree $3$, comparing degrees it must be that $v_2 = u_2$.
		However this then implies that $v_1$ and $v_2$ are adjacent to both $u$ and $v$, which contradicts that $v_1$ and $v_2$ have only one common neighbour.
		Therefore the vertices $u,v, u_1, u_2, v_1$ and $v_2$ are pairwise distinct.
		And so in $\Gamma-u$, $v$ has  degree two, and its neighbours are adjacent, have degree three and have only one common neighbour.
		Therefore $v$ is removable in $G-u$.
		
		To complete the proof it suffices to prove that any removable vertex in $\Gamma-u$ is removable in $\Gamma$.
		Let $w$ be a removable vertex in $\Gamma-u$.
		Since $\Gamma$ has maximum degree three, the neighbours of $w$, say $w_1$ and $w_2$ which have degree three in $\Gamma-u$, also have degree three in $\Gamma$.
		Therefore $w_1$ and $w_2$ are not adjacent to $u$
		
		Since no edges were added to $\Gamma$ to form $\Gamma-u$, if $w_1$ and $w_2$ have exactly one common neighbour in $\Gamma-u$ they must have one common neighbour in $\Gamma$.
		And since $w_1$ and $w_2$ are adjacent in $\Gamma - u$, by construction of $\Gamma-u$ they must also be adjacent in $\Gamma$.
		And so to prove $w$ is removable in $\Gamma$ it suffices to prove it has degree $2$ in $\Gamma$.
		
		Since $w$ is removable in $\Gamma-u$ it has degree two in $\Gamma-u$.
		The only way for it to have degree three in $\Gamma$ is if it is adjacent to $u$ in $\Gamma$.
		This would then imply, without loss of generality, that $w = u_1$.
		However in this case since $u_1$ and $u_2$ would be adjacent to both $u$ and $w_1$, contradicting that $u$ is removable in $\Gamma$.
		Therefore $w$ has degree two in $\Gamma$, which then implies it is removable in $\Gamma$. 
\end{proof}

Together Lemmas  \ref{lem:remove1} and \ref{lem:removePreserve} imply that removing the set of removable vertices of $\Gamma$ produces a reduced graph $\Gamma^R$ with the property that $\Gamma$ admits a partial orientation as a quasi-transitive mixed graph if and only if $\Gamma^R$ admits a partial orientation as a quasi-transitive mixed graph.

\begin{theorem}\label{thm:reduct3}
	Let $\Gamma$ be a graph with maximum degree three.
	Let $R$ be the set of removable vertices of $\Gamma$.
	The graph $\Gamma$ admits a partial orientation as a quasi-transitive mixed graph if and only if $\Gamma - R$ admits a partial orientation as a quasi-transitive mixed graph.
\end{theorem}

\begin{proof}	
The result follows directly by induction on $|R|$ by applying Lemmas \ref{lem:remove1} and \ref{lem:removePreserve}.
\end{proof}

\begin{figure}
	\includegraphics[scale=0.4]{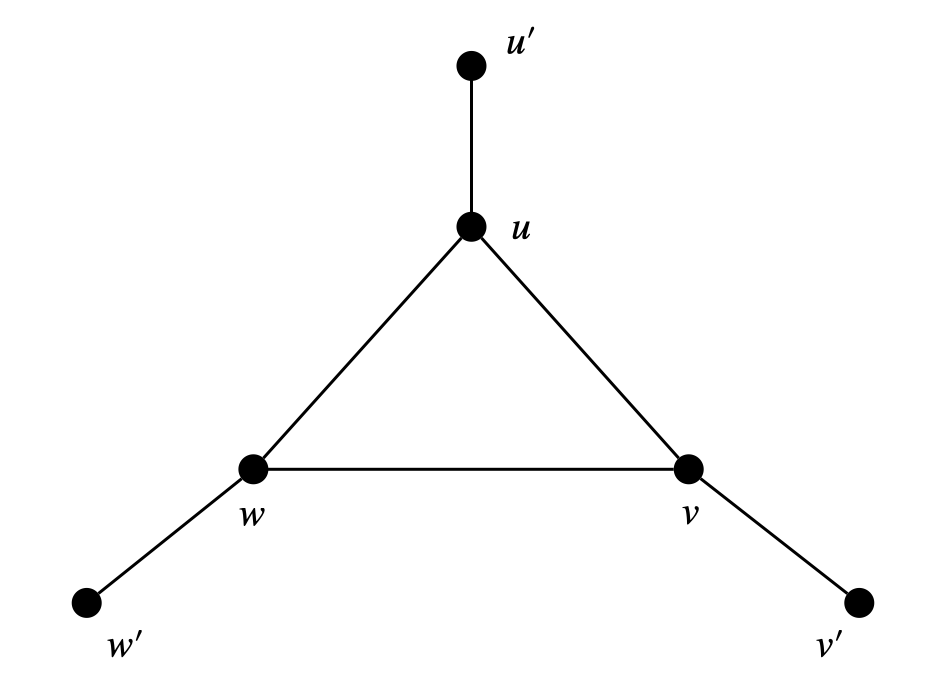}
	\caption{The graph $\Pi$, a forbidden subgraph for quasi-transitive mixed partial orientability of graphs with maximum degree three.}
	\label{fig:forbiddenCubic}
\end{figure}

We classify those reduced graphs with maximum degree three that arise as an oriented square by way of forbidden subgraphs.
Let $\Pi$ be the split graph on six vertices consisting of a copy of $K_3$ with a pendant edge on each vertex (see Figure \ref{fig:forbiddenCubic}).
Recall $\Pi$ is a forbidden subgraph for the family of graphs that admit a quasi-transitive orientation \cite{G88}.

\begin{lemma}\label{lem:noPi}
	Let $\Gamma$ be a graph with maximum degree three.
	If $\Pi$ is a subgraph of $\Gamma$, then $\Gamma$ does not admit a partial orientation as a quasi-transitive mixed graph.
\end{lemma}

\begin{proof}
	Let $\Gamma$ be a graph with maximum degree three containing $\Pi$.
	We proceed by contradiction.
	Let $G$ be a partial orientation of $\Gamma$ as a quasi-transitive mixed graph.
	
	Since $\Gamma$ has maximum degree three, none of $uu^\prime$, $vv^\prime$ and $ww^\prime$ are contained in a copy of $K_3$ in $\Gamma$.
	And so by Lemma \ref{lem:noK3}, each are oriented is an arc in $\Gamma$ and further each of $u$,$v$ and $w$ is a source or a sink.

	Without loss of generality, assume $u$ and $w$ are source vertices.
	Therefore $wu$ is an edge in $G$.
	Since $\Gamma$ has degree three, the only copy of $K_3$ in $\Gamma$ that contains $wu$ is the one with vertex set $\{u,v,w\}$.
	Therefore $uvw$ is a $2$-dipath (in some direction).
	Hence $v$ is not a source nor a sink, a contradiction.
\end{proof}

\begin{figure}
	\includegraphics[width=\linewidth]{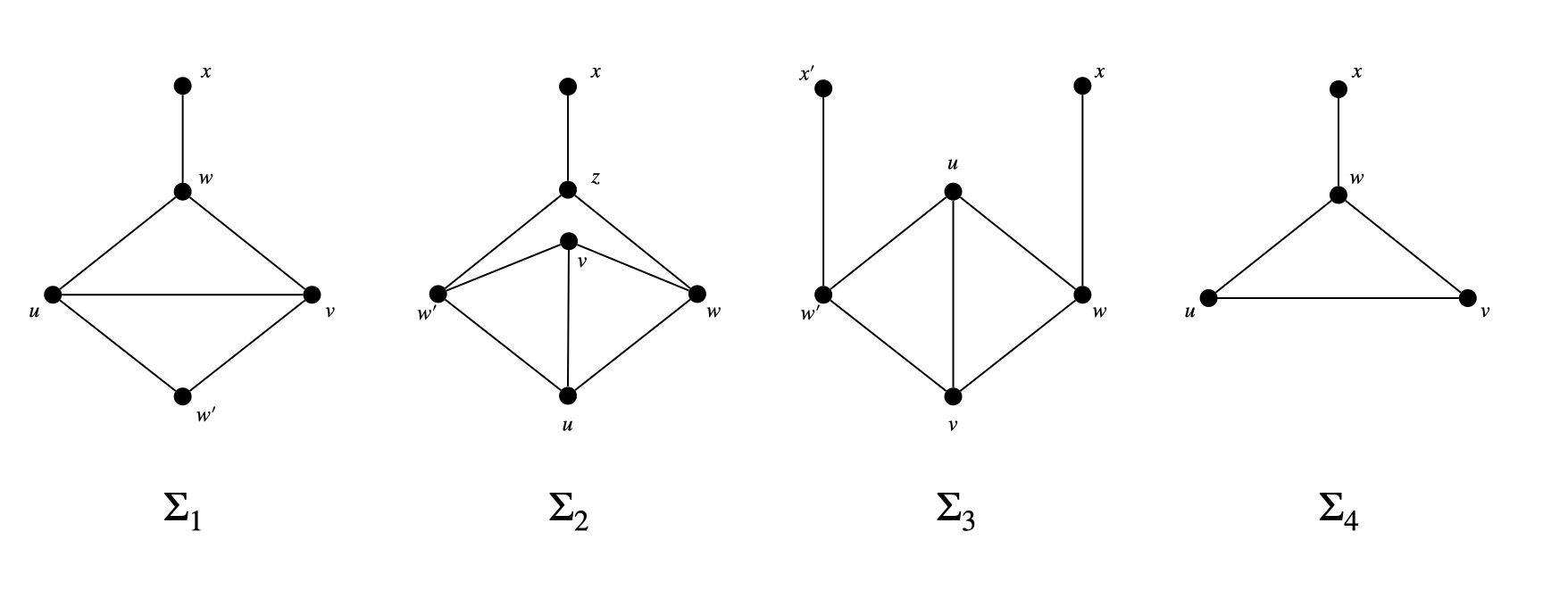}
	\caption{Configurations of copies of $K_3$ in graphs with maximum degree three that admit a partial orientation as a quasi-transitive mixed graph.}
	\label{fig:sigmasCubic}
\end{figure}

Lemma \ref{lem:noPi} restricts induced subgraphs containing of copies of $K_3$ in reduced graphs with maximum degree three that admit a partial orientation as a  quasi-transitive mixed graph.
The subsequent lemma characterises these subgraphs as those appearing in Figure \ref{fig:sigmasCubic}.

\begin{lemma}\label{lem:K3class}
Let $\Gamma$ be a reduced graph with maximum degree $3$ such that $\Gamma \neq K_3,K_4, K_4-e$
Consider a copy of $K_3$ in $\Gamma$ with vertices $u,v$ and $w$.
If $\Gamma$ is $\Pi$-free, then $u,v$ and $w$ are contained in an induced subgraph as shown in Figure \ref{fig:sigmasCubic}
\end{lemma}

\begin{proof}
	We proceed based on the degree of each of $u,v$ and $w$ in $\Gamma$
	
	If each vertex has degree three, then since $\Gamma$ is $\Pi$-free, without loss of generality, $u$ and $v$ have a second common neighbour, say $w^\prime$.
	Since $\deg(w) = 3$ and $\Gamma \neq K_4$, $w$ is not adjacent to $w^\prime$.
	If $\deg(w^\prime) = 2$, then $u,v$ and $w$ are configured as in $\Sigma_1$. 
	If $\deg(w^\prime) = 3$ and $ww^\prime \in E(\Gamma)$, then $u,v$ and $w$ are configured as in $\Sigma_2$.
	If $\deg(w^\prime) = 3$ and $ww^\prime \notin E(\Gamma)$, then $u,v$ and $w$ are configured as in $\Sigma_3$.
	
	If $u$ and $v$ have degree three and $w$ has degree two,  then $u,v$ and $w$ are configured as in $\Sigma_1$, exchanging the label $w$ for $w^\prime$, since $\Gamma \neq K_4-e$  and $\Gamma$ is reduced.
	
	Finally, consider the case where $u$ and $v$ have degree two. 
	Since $\Gamma \neq K_3$, $w$ must have degree three.
	And so $u,v$ and $w$ are configured as in $\Sigma_4$.
\end{proof}

\begin{figure}
	\includegraphics[width=\linewidth]{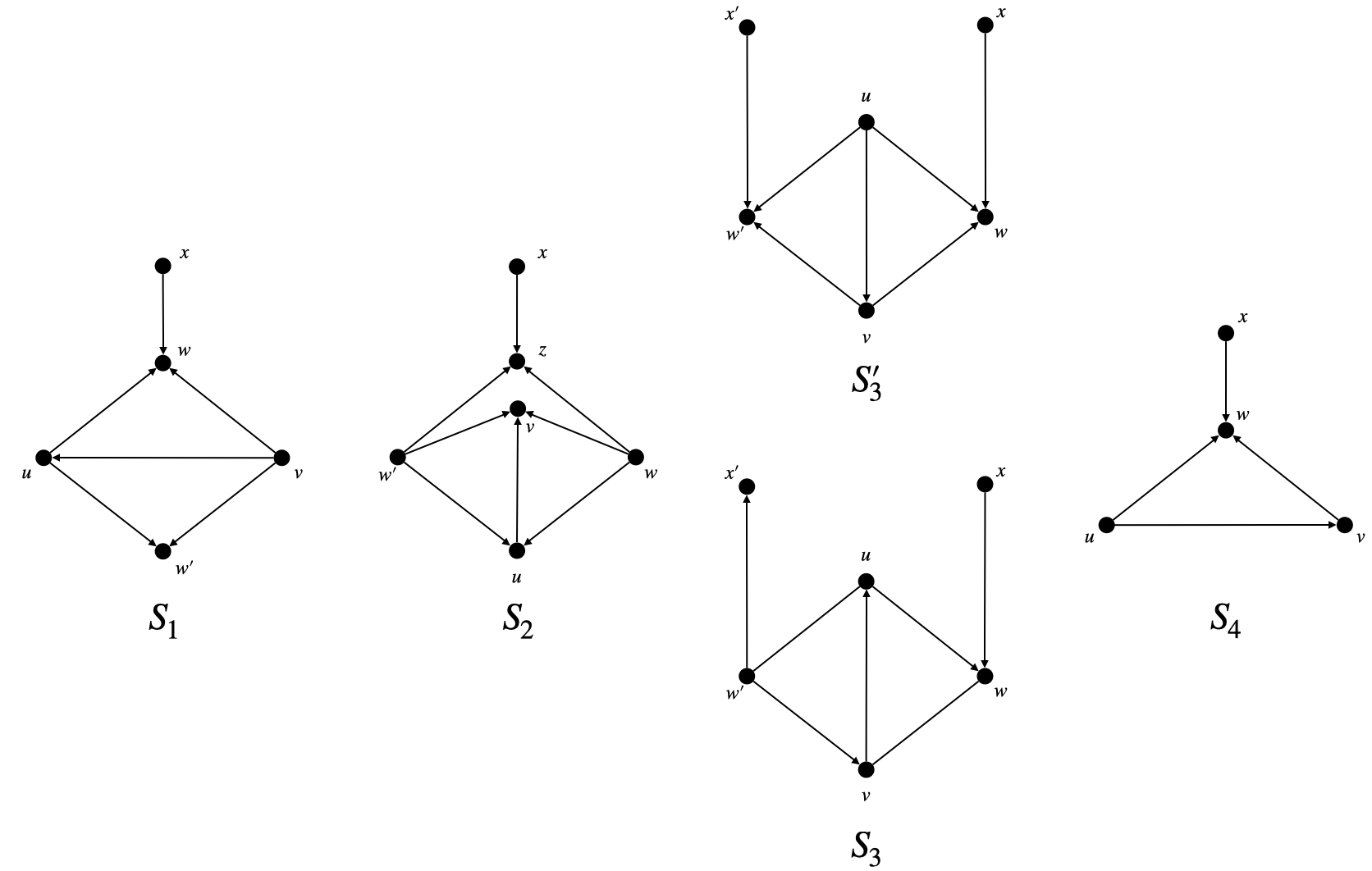}
	\caption{Partial orientations of graphs in Figure \ref{fig:sigmasCubic} as quasi-transitive mixed graphs.}
	\label{fig:essesCubic}
\end{figure}

Notice that if $\Sigma \in \{\Sigma_1,\Sigma_2,\Sigma_3,\Sigma_4\}$ is an induced subgraph of $\Gamma$, then $\Gamma$ contains an independent vertex cut.
Thus in determining if a graph containing a copy of $\Sigma_i$ admits a partial orientation as a quasi-transitive mixed graph we may appeal to the Vertex Cut Lemma.
With this strategy in mind, we note that Figure \ref{fig:essesCubic}, gives a partial orientation as a quasi-transitive mixed graph for each of the graphs in Figure \ref{fig:sigmasCubic}.

\begin{theorem}\label{thm:cubicMain}
	Let $\Gamma$ be a reduced graph with maximum degree three.
	The graph $\Gamma$ admits a partial orientation as a quasi-transitive mixed graph if and only if $\Gamma$ contains no copy of $\Pi$ and the subgraph formed by the set of edges of $\Gamma$ that are not contained in a copy of $K_3$ contains no odd cycle.
\end{theorem}

\begin{proof}
	Let $\Gamma$ be a reduced graph with maximum degree three.
	
	Assume $\Gamma$ admits a partial orientation as quasi-transitive mixed graph.
	By Lemma \ref{lem:noPi}, $\Gamma$ is $\Pi$-free.
	By Lemma \ref{lem:noK3}, the set of edges of $\Gamma$ that are not contained in a copy of $K_3$ is bipartite.
		
	To prove the converse we proceed by disproving the existence of a minimum counter example.
	Let $\Gamma$ be a reduced graph with maximum degree three that is $\Pi$-free so that the subgraph formed by the set of edges of $\Gamma$ that are not contained in a copy of $K_3$ is bipartite and $\Gamma$ does not admit a partial orientation as a quasi-transitive mixed graph.
	Among all choices for $\Gamma$,  choose  one with the fewest number of vertices.
	
	Notice that each of $K_3,K_4$ and $\Gamma = K_4-e$ are comparability graphs, and thus admit an partial orientation as a quasi-transitive mixed graph.
	Therefore $\Gamma \neq K_3,K_4, K_4-e$.
	
	Let $Y$ be the set of edges contained in no copy of $K_3$ in $\Gamma$.
	Notice $Y \neq \emptyset$, as otherwise by Theorem \ref{thm:bigGirth} $\Gamma$ admits a partial orientation as a quasi-transitive mixed graph.
	Since $Y \neq \emptyset$, by Lemma \ref{lem:K3class} $\Gamma$ contains one of $\Sigma_1, \Sigma_2, \Sigma_3$ or $\Sigma_4$ as an induced subgraph.
	
	Let $\Sigma^\prime$ be an induced subgraph of $\Gamma$ such that $\Sigma^\prime \in \{\Sigma_1, \Sigma_2, \Sigma_3, \Sigma_4\}$.
	Consider the graph $\Gamma - X$ where 
	\begin{itemize}
		\item $X = \{u,v,w^\prime\}$ if $\Sigma^\prime = \Sigma_1$;
		\item $X = \{u,v,w,w^\prime\}$  if $\Sigma^\prime = \Sigma_2$; and
		\item $X = \{u,v\}$ if $\Sigma^\prime = \Sigma_3,\Sigma_4$.
	\end{itemize}

	By minimality of $\Gamma$ and Lemma \ref{lem:removePreserve}, $\Gamma-X$ admits a partial orientation $G$ as a quasi-transitive mixed graph.
	
We extend this partial orientation of  $\Gamma-X$ to one of $\Gamma$ as follows:
\begin{itemize}
	\item if $\Sigma^\prime = \Sigma_1$, then orient the edges $uw, wv, vw^\prime, w^\prime u$ as in $S_1$ in Figure \ref{fig:essesCubic} if $xw \in A(G)$ or the converse if $wx \in A(G)$;
	\item if $\Sigma^\prime = \Sigma_2$, then orient the edges $uw, wz, zw^\prime, w^\prime u, uv, vw$ and $vw^\prime$ as in $S_2$ in Figure \ref{fig:essesCubic} if $xz \in A(G)$ or the converse if $zx \in A(G)$
	\item if $\Sigma^\prime = \Sigma_3$, then orient the edges $w^\prime v, vu$ and $uv$ as in $S_3$ in  Figure \ref{fig:essesCubic} if $w^\prime x^\prime, xw \in A(G)$ or the converse if $x^\prime w^prime, xw \in A(G)$
	\item if $\Sigma^\prime = \Sigma_3$, then orient the edges $w^\prime v, vw,wu,uw$ and $uv$ as in $S_3^\prime$ in Figure  \ref{fig:essesCubic} if $x^\prime w^\prime, xw \in A(G)$ or the converse if $w^\prime x^prime, wx \in A(G)$
	\item if $\Sigma^\prime = \Sigma_4$, then orient the edges $uv,vw$ and $wu$ as in $S_4$ in Figure \ref{fig:essesCubic} if $xw \in A(G)$ or the converse if $wx \in A(G)$
\end{itemize}
	This construction contradicts that $\Gamma$ admits no partial orientation as a quasi-transitive mixed graph, which completes the proof.
	\end{proof}

By noticing that a cubic graph is necessarily reduced, we arrive at the following classification of cubic graphs that admit a partial orientation as a quasi-transitive mixed graph

\begin{corollary}
	A cubic graph $\Gamma$ admits a partial orientation as a quasi-transitive mixed graph if and only if $\Gamma$ is $\Pi$-free that the subgraph of $\Gamma$ formed by the set of edges of $\Gamma$ that are not contained in a copy of $K_3$ contains no odd cycle.
\end{corollary}

As the collection of edges that are contained in no copy of $K_3$ and removable vertices can be found in polynomial time, Theorem \ref{thm:cubicMain} implies that it is Polynomial to decide if a graph with maximum degree three is an oriented graph square.


\section{Decision Problems and Computational Complexity}\label{sec:complexity}

We prove the problem of deciding if a graph is an oriented graph square is NP-complete via reduction from the problem of monotone not-all-equal $3$-satisfiability.

\underline{Monotone NAE3SAT}\\
\emph{Instance}: A monotone boolean formula $Y = (L,C)$ in conjunctive normal form with three variables in each clause\\
\emph{Question}: Does there exist a not-all-equal satisfying assignment for the elements of $L$?

\begin{theorem}\cite{S78}
	Monotone NAE3SAT is NP-complete.
\end{theorem}

\underline{MIXEDQT}\\
\emph{Instance}: A graph $\Gamma$\\
\emph{Question} Does $\Gamma$ admit a partial orientation as a mixed quasi-transitive graph?

\underline{UNDIRSQUARE}\\
\emph{Instance}: A graph $\Gamma$\\
\emph{Question} Does $\Gamma$ arise as an undirected square of an oriented graph?

Our proof proceeds by constructing a gadget graph for each clause so that whether a particular vertex in the gadget is a source or a sink corresponds to the truth value of a literal in the clause.

\begin{figure}
	\includegraphics[width=0.5\linewidth]{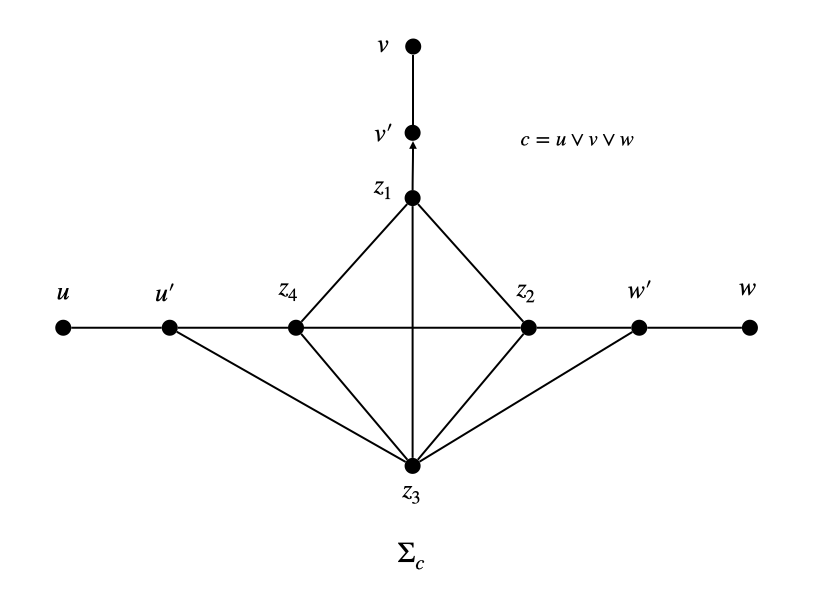}
	\caption{The Clause Graph}
	\label{fig:clauseGraph}
\end{figure}

Let $Y = (L,C)$ be an instance of Monotone NAE3SAT.
Let $c \in C$ and let $\Sigma_c$ be the graph shown in Figure \ref{fig:clauseGraph}.
We refer to $\Sigma_c$ as  a \emph{clause graph}.

\begin{figure}
	\includegraphics[width = \linewidth]{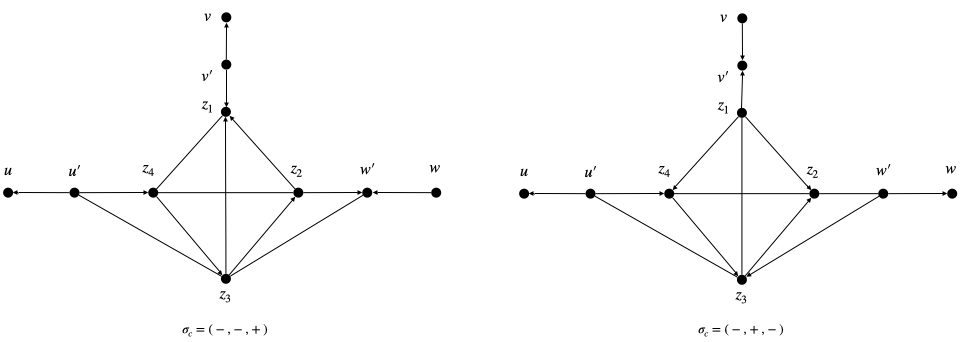}
	\caption{Partial orientations of a clause graph with signatures $\sigma_c = (-,-,+)$ and $\sigma_c = (-,+,-)$}
	\label{fig:clauseSignature}
\end{figure}

Let $S_c$ be a partial orientation of $\Sigma_c$ in which the edges $uu^\prime, vv^\prime, ww^\prime$ are all oriented as arcs in some direction.
We associate with $S_c$ a vector $\sigma_{c} = (\sigma_{u,c},\sigma_{v,c},\sigma_{w,c}) \in \{-,+\}^3$ where $\sigma_{u,c}$ (respectively $\sigma_{v,c}, \sigma_{w,c}$) is $+$ when $u$ (respectively $v, w$) is a source in $S_c$ and $-$ when $u$ (respectively $v,w$) is a sink in $S_c$.
We call $\sigma_c$ the \emph{signature} of $S_c$.
See Figure \ref{fig:clauseSignature} for an example.
We will use the signature of a partial orientation of a clause graph to denote the truth values of the variables in the clause graph.
The literal $u$ (respectively $v$ and $w$) will be TRUE in $C$ if and only if $\sigma_{u,c}$ (respectively, $\sigma_{v,c}$ and $\sigma_{w,c}$) is $+$ .

Previewing our approach to an NP-completeness proof, we take a moment to justify collapsing discussion of not-all-equal cases to one where we may restrict our consideration to satisfying assignments of a YES instance of not-all-equal satisfiability where we may assume without loss of generality that for $c = u \vee v \vee w$ we have that $u$ and $v$ are FALSE and $w$ is TRUE, or $u$ and $w$ are FALSE and $v$ is TRUE.

Consider the two partial orientations given in Figure \ref{fig:clauseSignature}.
Reversing the orientation of each of the arcs yields partial orientations with  $\sigma_c = (+,+,-)$ and $\sigma_c = (+,-,+)$.
Moreover, notice the existence of an automorphism of $\Sigma$ that swaps $u$ and $w$.
Applying this automorphism to the partial orientation on the left permits us to generate partial orientations with $\sigma_c = (-,+,+)$ and $\sigma_c = (+,-,-)$.
Thus there are partial orientations of $\Sigma$ as a quasi-transitive mixed graphs for each signature in $\{-,+\}^3 \setminus \{(-,-,-), (+,+,+)\}$.
The signature of a clause graph for $c = u \vee v \vee w$ will give the truth value of the literals $u$, $v$ and $w$ in $c$.

\begin{lemma}\label{lem:orientSigma}
	There exists a partial orientation $S_c$ of $\Sigma_c$ as a quasi-transitive mixed graph if and only if  $\sigma_c \in \{-,+\}^3 \setminus \{(-,-,-), (+,+,+)\}$	
\end{lemma}

\begin{proof}
	Consider  $\sigma_c \in \{-,+\}^3 \setminus \{(-,-,-), (+,+,+)\}$.
	Figure \ref{fig:clauseSignature} gives partial orientations $S_c$ of $\Sigma_c$ as a quasi-transitive mixed graph for $\sigma_c \in \{ (-,-,+), (-,+,-)\}$.
	The other cases are obtained by applying the automorphism of $\Sigma$ that swaps $u$ and $w$ and/or reversing the orientation of every arc.
	
	To complete the proof it suffices to prove there is no partial orientation  of $\Sigma_c$ as a quasi-transitive mixed graph with $\sigma_c  = (-,-,-)$ or $\sigma_c  = (+,+,+)$.
	By computer search (see code in Appendix \ref{app:code}) no such partial orientation exists.	
\end{proof}

To prove MIXEDQT is NP-hard, given an instance of  $Y = (L,C)$ of Monotone NAE3SAT we construct a graph, $\Gamma_Y$, that admits a partial orientation as a  quasi-transitive mixed graph if and only if $Y$ is a YES instance of Monotone NAE3SAT.

We construct $\Gamma_Y$ from the disjoint union of the following graphs:
\begin{itemize}
\item for each $c \in C$, a copy of $\Sigma_c$; and
\item for each $x \in L$, a path $P_x$ with $2|C|+2$ vertices: $x_0,x_1,\dots, x_{2|C| +2}$
\end{itemize}
by identifying vertices as follows.
 
\begin{figure}
	\includegraphics[scale=0.6]{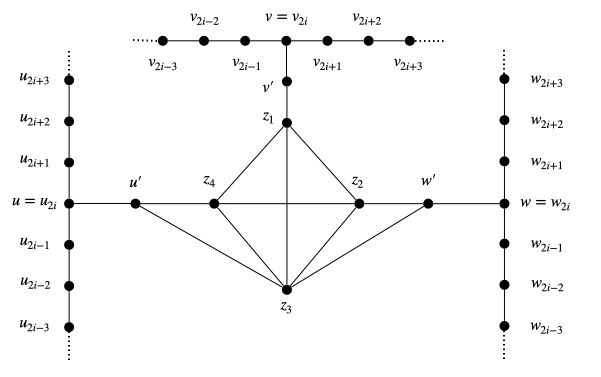}
	\caption{A clause subgraph in the construction of $\Gamma_Y$ for  the clause $c_i = u \vee v \vee w$}
	\label{fig:gammaConstruct}
\end{figure}

For every $1 \leq k \leq |C|$ with $c_k = u \vee v \vee w$, identify vertex $u_{2k}$ (respectively $v_{2k}$ and $w_{2k}$) in $P_u$ (respectively $P_v$ and $P_w$) with vertex $u$ (respectively $v$ and $w$) in $\Sigma_{c_k}$. (See Figure \ref{fig:gammaConstruct}).
In other words,  if $x\in L$ appears in clause $c_k \in C$, then vertex $x$ in $\Sigma_k$ is identified with vertex $x_{2k}$ in $P_x$.

Notice that for every literal $x \in L$, no edge in $P_x$ is contained in a copy of $K_3$ in $\Gamma_Y$.
Therefore  by Lemma \ref{lem:noK3}, in every partial orientation of $\Gamma_Y$ as a quasi-transitive mixed graph, every such edge is oriented as an arc.
And so in every every partial orientation of $\Gamma_Y$ as a quasi-transitive mixed graph, every path $P_x$ is oriented as alternating path where every vertex of the form $x_{2k}$ is a source (respectively, a sink) if and only if $x_0$ is a source (respectively, a sink).

\begin{lemma}\label{lem:altPath}
	Let $Y = (L,C)$ be an instance of Monotone NAE3SAT and consider $x \in L$.
	In every partial orientation of $\Gamma_Y$ as a quasi-transitive mixed graph
	\begin{enumerate}
		\item the path $P_x$ is oriented as an alternating path; and
		\item for every $1 \leq k \leq |C|/2$, $x_{2k}$ is a source if and only if $x_0$ is a source.
	\end{enumerate}
\end{lemma}

\begin{proof}
	Let $Y = (L,C)$ be an instance of Monotone NAE3SAT and let $G$ be a partial orientation of $\Gamma_Y$ as a quasi-transitive mixed graph.
	Consider $x \in L$.
	By construction, no vertex of $P_x$ is contained in a copy of $K_3$.
	And so by Lemma \ref{lem:noK3}, every edge of this path is oriented as an arc in $G$ and every vertex of this path is either a source or a sink.
	Noting that there cannot be an arc with both ends at a  source or sink vertex, it then follows that $P_x$ is oriented as an alternating path in $P_x$.
	From this it follows direct that for every $1 \leq k \leq |C|/2$, $x_{2k}$ is a source if and only if $x_0$ is a source.
\end{proof}

In our mapping between YES instances of Monotone NAE3SAT and partial orientations of $\Gamma_Y$ as a quasi-transitive mixed graph, $x_0$ as a source (respectively, a sink) corresponds to $x$ being TRUE (respectively, FALSE) in a satisfying assignment of the instance of Monotone NAE3SAT.

\begin{lemma}\label{lem:mainComplex}
	An instance $Y = (L,C)$ of Monotone NAE3SAT is a YES instance if and only if $\Gamma_Y$ admits a partial orientation as a mixed quasi-transitive graph.
\end{lemma}

\begin{proof}
	Let $Y = (L,C)$ be a YES instance of Monotone NAE3SAT.
	Consider the following partial orientation  of $G$ of $\Gamma_Y$:
	\begin{itemize}
		\item for each $x \in L$, if $x$ is true, then orient $P_x$ as an alternating path where $x_0$ is a source;
		\item for each $x \in L$, if $x$ is false, then orient $P_x$ as an alternating path where $x_0$ is a sink; and
		\item for each clause $c = u \vee v \vee w$, partially orient $\Sigma_c$ as a quasi-transitive mixed oriented graph whose signature matches the truth values of the literals $u, v$ and $w$.
	\end{itemize}
	Note that by Lemmas \ref{lem:orientSigma} and  \ref{lem:altPath} such an partial orientation necessarily exists.
	(See Figure \ref{fig:GConstruct} for the case $c_i = u \vee v \vee w$ where $u$, $v$ and $w$ are respectively FALSE, FALSE and TRUE.) 
	
	\begin{figure}
		\includegraphics[scale=0.6]{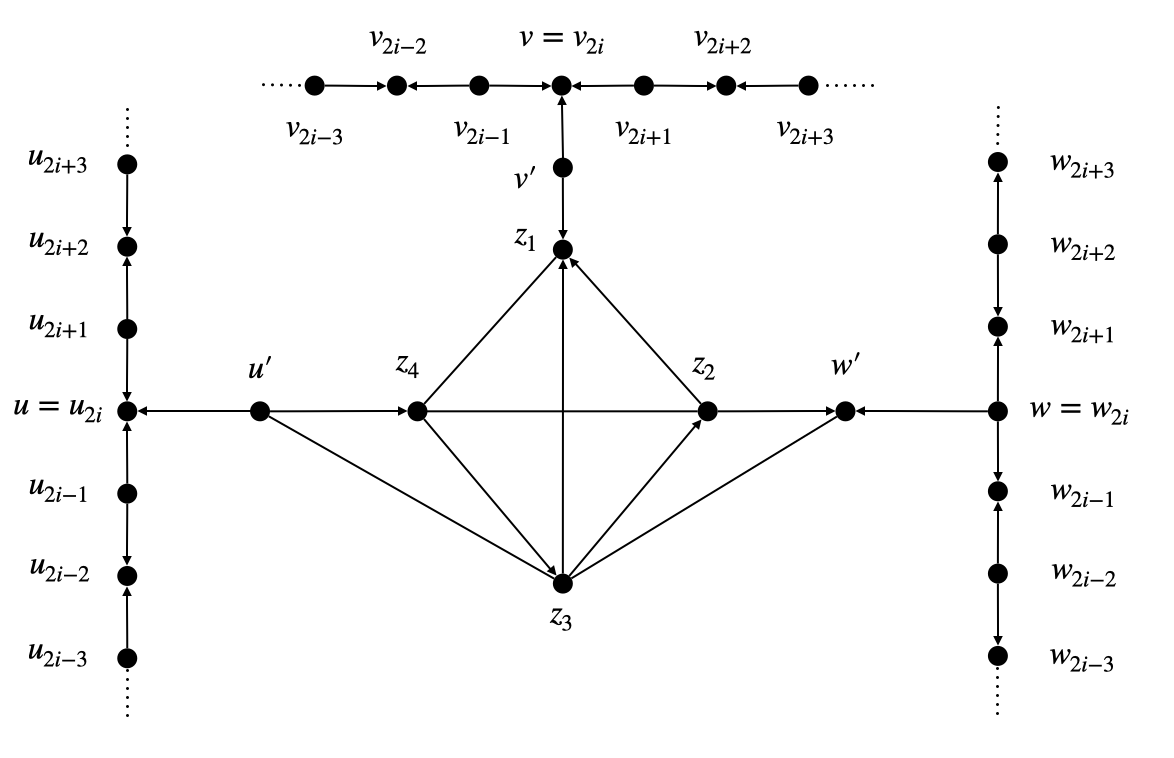}
		\caption{A clause subgraph in the partial orientation of $G$ in the proof of Theorem \ref{lem:mainComplex} where  $u$, $v$ and $w$ are respectively FALSE, FALSE and TRUE.}
		\label{fig:GConstruct}
	\end{figure}
	
	We claim $G$ is a quasi-transitive mixed graph.
	We proceed by appealing to the definition of quasi-transitive mixed graph.
	
	Let $y \in V(\Gamma)$.
	We prove $y$ is not the centre vertex of any induced $2$-dipath.
	We proceed in cases based on whether $y$ is contained in a clause graph.
	
	If there exists $c \in C$ such that $y \in V(\Sigma_c)$, then, without loss of generality, $y$ is configured as one of the vertices in Figure \ref{fig:clauseSignature}.
	If $y \neq u,v,w$, then by observing the partial orientations in vertices in Figure \ref{fig:clauseSignature}, we see that $y$ is not the centre vertex of an induced $2$-dipath.
	If $y \in \{u,v,w\}$, then by construction of $\Gamma_Y$ there exists $1 \leq i \leq |C|/2$ such that $y = u_{2i}$, $y = v_{2i}$ or $y = w_{2i}$.
	By the Vertex Cut Lemma with $I = \{u_{2i}, v_{2i},w_{2i}\}$, $y$ is a source or a sink.
	Therefore $y$ is not the centre vertex of an induced $2$-dipath in $G$.
	
	If $y$ is contained in no clause graph, then there exists $x \in L$ such that $y = x_0$ or $y = x_{2i+1}$ for $1 \leq i \leq |C|/2$.
	In either case, we note that by Lemma \ref{lem:altPath}, $y$ is a source or a sink in $G$.
	And so we conclude $y$ is not the centre vertex of an induced $2$-dipath in $G$.
		
	Consider $e \in E(G)$.
	By construction, there exists $c \in C$ such that $e \in E(\Sigma_c)$.
	Without loss of generality, assume $\sigma_c \in \{ (-,-,+), (-,+,-)\}$.
	Therefore $\Sigma_c$ is partially oriented as one of the mixed graphs in Figure \ref{fig:clauseSignature}.
	We observe that in either case the end points of every edge are also the endpoints of a $2$-dipath. 
	Therefore $G$ is a quasi-transitive mixed graph.
		
	Let $G$ be a partial orientation of $\Gamma$ as a quasi-transitive mixed graph.
	We assign the truth value of the literal $x$ based on whether $x_0 \in V(G)$ is a source or sink.
	Let $f:L \to \{TRUE, FALSE\}$ so that $f(x) = TRUE$ if and only if $x_0$ is a source in $G$.
	Notice that since $x_0$ has degree $1$ in $\Gamma$, $x_0$ is necessarily a source or a sink in $G$.
	And so $f(x) = FALSE$ if and only if $x_0$ is a sink in $G$.
	
	Notice that for any clause $c_i = u \vee v \vee w$, the set $I = \{u_{2i}, v_{2i}, w_{2i}\}$ satisfies the hypothesis of the Vertex Cut Lemma.
	Therefore for each $x \in L$ and each $1 \leq i \leq |C|/2$, $x_{2i}$ is a source or a sink.
	(Note that if $x$ is not contained in clause $c_i$, then $x_{2i}$ is a source or a sink by Lemma \ref{lem:altPath}.
	By Lemma \ref{lem:altPath}, $x_{2i}$ is a source if and only if $x_0$ is a source.
	Therefore for every $c \in C$ with $c = u \vee v \vee w$, we have $\sigma_{c,u} = +$ if and only if $u_0$ (respectively $v_0$ and $w_0$) is a source.
	By Lemma \ref{lem:orientSigma},  $\sigma_c \notin \{(-,-,-),(+,+,+)\}$.
	Therefore $f$ is a satisfying assignment for $Y$.	
\end{proof}

\begin{theorem}\label{thm:mainComplex}
	The decision problem MIXEDQT is NP-complete.  
	The problem is Polynomial when restricted to inputs of maximum degree $3$, but remains NP-compete when restricted to inputs of maximum degree $k$ for all $k \geq 5$.
\end{theorem}

\begin{proof}
	The reduction is from Monotone NAE-3SAT, noting that given an instance $Y$ of Monotone NAE3SAT, $\Gamma_Y$ can be constructed in polynomial time.
	By Lemma \ref{lem:mainComplex}, $\Gamma_Y$ admits a partial orientation as a quasi-transitive mixed graph if and only if $Y$ is a YES instance.
	Therefore the problem is NP-complete.
	 
	For every instance $Y$,  $\Gamma_Y$ has maximum degree five, thus problem remains NP-complete when restricted to inputs with maximum degree five.
	Further, $\Gamma_Y$ can be modified by adding arbitrarily many pendants to a vertex $x_0$ without changing the correspondence between YES instances and partial orientations as quasi-transitive mixed graphs.
	Therefore the problem remains NP-complete when restricted to inputs with maximum degree $k$ for any $k \geq 5$.
	By Theorem \ref{thm:cubicMain}, the problem is Polynomial when restricted to inputs with maximum degree three.
\end{proof}

\begin{corollary}
	The decision problem UNDIRSQUARE is NP-complete. 
	The problem is Polynomial when restricted to inputs of maximum degree $3$, but remains NP-compete when restricted to inputs of maximum degree $k$ for all $k \geq 5$.
\end{corollary}

We notice here the gap in the analysis for the case $\Delta = 4$.
It is possible that the analysis done in Section \ref{sec:cubic} can be undertaken for graphs with maximum degree four.
However in this case $\Pi$ is not a forbidden subgraph, as demonstrated by the example in Figure \ref{fig:exampleFig}.

Let $\mathcal{F}$ be a family of oriented graphs. 
Consider the following decision problem.

\underline{$\mathcal{F}$-UNDIRSQUARE}\\
\emph{Instance}: A graph $\Gamma$.\\
\emph{Question} Does $\Gamma$ arise as an undirected square of an oriented graph $G \in \mathcal{F}$?

Using the proof of Theorem \ref{thm:mainComplex}, we consider $\mathcal{F}_k$-UNDIRSQUARE where $\mathcal{F}_k$ is the family of orientations of graphs with maximum degree $k$.

\begin{theorem}\label{thm:UnDirMaxk}
	Let $k$ be a fixed positive integer and let $\mathcal{F}_k$ be the family of orientations of graphs with $\Delta = k$.
	For $k \geq 3$ the decision problem $\mathcal{F}_k$-UNDIRSQUARE is NP-complete.
\end{theorem}

\begin{proof}	
	The reduction is from Monotone NAE3SAT, noting that given an instance $Y$ of Monotone NAE3SAT, $\Gamma_Y$ can be constructed in polynomial time.
	
	Notice that in any YES instance $Y$ of Monotone NAE-3SAT, there is a corresponding partial orientation of $\Gamma_Y$ as a quasi-transitive mixed graph is one in which every vertex is incident with at most three arcs (see Figure \ref{fig:essesCubic}).
	Therefore $Y$ is a YES instance if and only if there exists an oriented graph $\ovr{H}$ with maximum degree three such that $\Gamma_Y = U(\ovr{H}^2)$.
	And so it is NP-complete to decide of a graph is an oriented graph square of an oriented graph with maximum degree three.
	
	For fixed $k >3$, as in the proof of Theorem \ref{thm:mainComplex}, we can modify the construction of $\Gamma_Y$ by adding arbitrarily many pendants to a vertex $x_0$ without changing the correspondence between YES instances and partial orientations as quasi-transitive mixed graphs.
	In these partial orientations every edge incident with $x_0$ is oriented as an arc.
	Therefore  $Y$ is a YES instance if and only if there exists an oriented graph $\ovr{H}$ with maximum degree $k$ such that $\Gamma_Y = U(\ovr{H}^2)$.
	And so it is NP-complete to decide of a graph is an oriented graph square of an oriented graph with maximum degree $k$ for any fixed $k \geq 3$.
\end{proof}

\section{Discussion and Future Work}\label{sec:discussion}

As with other extensions of graph theoretic concepts to the study of oriented graphs, the landscape of open areas for investigation on this topic is vast.
For example, one may consider graphs arising from orientations of trees.
As with the case for graphs, we expect there to be a classification that yields a polynomial time algorithm for identifying graphs that arise as undirected squares of orientations of trees.

Recall that one may extend the definition of graph colouring to oriented graphs by way of oriented graph homomorphism (see \cite{S16}.
In this context, an oriented $k$-colouring of an oriented graph $\ovr{G}$ is a homomorphism to a tournament on $k$ vertices.
The oriented chromatic number, denoted $\chi_o(\ovr{G})$ is the least integer $k$ such that $\ovr{G}$ admits an oriented $k$ colouring.
In every oriented colouring vertices at directed distance $2$ necessarily are assigned different colours.
Hence $\chi(U(\ovr{G}^2)) \leq \chi_o(\ovr{G})$.
And so the oriented chromatic number of the underlying oriented graph gives an upper bound for the chromatic number of undirected square.

For many fundamental families of oriented graph, the problem of computing the oriented chromatic number remains open.
For example, the best known upper bound on the chromatic number of the family of orientations of planar graphs is $80$ \cite{RASO94}.
This bound is not expected to be tight; there are no known constructions for oriented planar graphs for which the chromatic number exceeds $16$ \cite{BDS17}.
The example of an orientation of a planar graph $\ovr{G}$ with oriented chromatic number $16$ in fact has $16$ vertices and thus satisfies $U(\ovr{G}^2) = K_{16}$. 
Hence $\chi(U(\ovr{G}^2)) = 16 =  \chi_o(\ovr{G})$.

Oriented cliques (i.e., oriented graphs for which $\chi_o(G) = |V(G)|$) present an example of pairs of non-isomorphic oriented graphs that have the same undirected square.
For example, $K_5$ is an oriented graph square of both the directed five cycle and the transitive tournament on $5$ vertices.
And so $K_5$ does not have a unique square root with respect to oriented raphs
Graphs with a unique square root (with respect to graphs) are known to exist
For example, tree squares are unique and, more generally, as are squares of graphs with girth at least $6$ \cite{R60,A09}.
The former of these results does not extend to graphs arising as undirected squares of orientations of trees.
The directed path on three vertices and an orientation of $K_{1,3}$ in which the vertex of degree three is not a source or a sink generate the same undirected square.

\bibliographystyle{abbrv}
\bibliography{references}

\appendix

\section{Code for Proof of Lemma \ref{lem:orientSigma}} \label{app:code}

This code checks that $\Sigma$ has no partial orientation as a quasi-transitive graph with $\sigma_c = (+,+,+)$ or $\sigma=(-,-,-)$. We do this by generating all partial orientations as quasi-transitive graph via brute force search of all mixed graphs that can be formed from $\Sigma$.

We generate mixed graphs as follows: for every subset of the edge set of $\Sigma$ we assign these edges to remain edges in the mixed graphs and the iterate over all orientations of the remaining edges, testing for each orientation if the constructed mixed graph is quasi-transitive. 
Those that are valid are stored. 
We then iterate over all of the valid partial orientations and check that none has $\sigma_c = (+,+,+)$ or $\sigma_c = (-,-,-$).
The expected output is \verb|0,0|

\begin{python}

sigma = graphs.PathGraph(6)
D = Graph(9)
G = Graph(9)
dict=[] #This list stores all possible partial orientations of \Sigma as a mQT
ttt=[] #This list stores all possible partial orientations of \Sigma as a mQT with \sigma_c = (+,+,+)
fff=[] #This list stores all possible partial orientations of \Sigma as a mQT with \sigma_c = (-,-,-)

sigma.add_edges([(2,7),(7,3),(6,7),(1,6),(6,4),(2,6),(3,6),(7,8)])

#generate all subsets of edges of \Sigma.     
edgeSet = sigma.edges()
subgraphs = [[]]
for e in edgeSet:
	for i in range(len(subgraphs)):
		subgraphs += [subgraphs[i]+[e]]

#Iterate over all choice of edges
for subEdge in subgraphs:
	D = Graph(9)
	G = Graph(9)
	for  e in sigma.edges():
		G.add_edge(e)
	D.add_edges(subEdge)
	G.delete_edges(subEdge)

		#Iterate over all orientations of the edges that were not chosen
		allOrientations = D.orientations()
		for  J in allOrientations:
			isValid = isMQT(G,J)
			if isValid == true:
				dict.append((G,J))

# At this stage, dict holds all mQTs that can be formed from \Sigma. 
# We now count the number that have \sigma_c = (+,+,+) or \sigma_c = (-,-,-).

for i in range(len(dict)):
	D = dict[i]
	if (0,1,None) in dict[i][1].edges() and (7,8,None) in dict[i][1].edges() and (5,4,None) in dict[i][1].edges():
		ttt.append(dict[i][1])
	elif (1,0,None) in dict[i][1].edges() and (8,7,None) in dict[i][1].edges() and (4,5,None) in dict[i][1].edges():
		fff.append(dict[i][1])

print(len(ttt),len(fff))
\end{python}

\end{document}